\documentclass[11pt,oneside,a4paper]{article}
\usepackage{syntonly}
\usepackage{geometry}
\usepackage{mathrsfs}
\usepackage{makeidx}
\makeindex
\usepackage{amsmath}
\usepackage{amssymb}
\usepackage{amsmath}
\usepackage{amsfonts}
\usepackage{amssymb}
\usepackage{esvect}
\usepackage{algorithm}
\usepackage{bm}
\usepackage{tikz}

\usepackage{cleveref}

\usepackage{graphicx}
\usepackage{epstopdf}

\DeclareMathOperator{\arctanh}{arctanh}

\newtheorem{theorem}{Theorem}

\newcommand{\qed}{\nobreak \ifvmode \relax \else
      \ifdim\lastskip<1.5em \hskip-\lastskip
      \hskip1.5em plus0em minus0.5em \fi \nobreak
      \vrule height0.30em width0.4em depth0.25em\fi}

 \author{Safari Mukeru\\
\footnotesize{\em Department of Decision Sciences}\\ 
\footnotesize{University of South Africa, P. O. Box 392, Pretoria, 0003. South Africa}\\
\footnotesize{e-mail:mukers@unisa.ac.za}}

\title{{Average number of real zeros of random  algebraic polynomials defined by the increments of fractional Brownian motion}}

\date{}

\begin{document}

\maketitle

\pagenumbering{arabic}

\begin{abstract}
The study of random polynomials  has a long and  rich history.  This paper  studies random algebraic polynomials $P_n(x) = a_0 + a_1 x + \ldots + a_{n-1} x^{n-1}$ where the coefficients $(a_k)$ are correlated random variables taken as the increments $X(k+1) - X(k)$, $k\in \mathbb{N}$, of a fractional Brownian motion $X$ of Hurst index $0< H < 1$. This reduces to the classical setting of independent coefficients for $H = 1/2$. We obtain that the average number of the real zeros of $P_n(x)$ is~$\sim K_H \log n$, for large $n$, where $K_H = (1 + 2 \sqrt{H(1-H)})/\pi$ (a generalisation of a classical result  obtained by Kac in 1943).  Unexpectedly, the parameter $H$ affects only the number of positive zeros, and the number of real zeros of the polynomials corresponding to fractional Brownian motions of indexes $H$ and $1-H$ are essentially the same. The limit case $H = 0$ presents some particularities: the average number of positive zeros converges to a constant.  These results shed some light on the nature of fractional Brownian motion on the one hand and on the behaviour of real zeros of random polynomials of dependent coefficients on the other hand.  
\end{abstract}


\section{Introduction}
Given a sequence $\{a_0, a_1, a_2, \ldots\}$ of independent, identically distributed random variables, consider the random polynomial
              $$P_n(x) = a_0 + a_1 x + \ldots+ a_{n-1} x^{n-1}.$$
The study of random polynomials $P_n(x)$ is a well-known problem. We shall denote  by $N_n$ the number of real zeros of $P_n(x)$ and $E_n = \mathbb{E}(N_n)$ its average.  In 1943, Kac \cite{Kac_1943} proved that if the random variables $a_i$ have the standard normal distribution, then 
      \begin{eqnarray} \label{kacform2}
  E_n \sim (2/\pi) \log n,\,\, n \to \infty.
\end{eqnarray}
      Later Kac \cite{Kac_1949}   
extended this asymptotic estimate to the case where  the random variable $a_i$ is uniformly distributed on $[-1, 1]$ and 
Erd\"os and Offord \cite{Erdos_Offord} to the case where $a_i$ is uniform on $\{-1, 1\}$. We refer the reader to Bharucha-Ried and Sambandham \cite{Bharucha-Reid_Sambandham}),  Dembo et al. \cite{Dembo_et_al}, Ibragimov and Maslova \cite{Ibragimov_Maslova} and  Ibragimov and Zeitouni \cite{Ibragimov_Zeitouni} for historical developments and further results on the subject. 
Recently Matayoshi \cite{Matayoshi} considered the case where the coefficients $a_0, a_1, \ldots$ are dependent but form a stationary sequence of standard normal distributions  and obtained essentially that $E_n \sim \frac{2}{\pi} \log n$ for $n \to \infty$.  Nezakati and  Farahmand \cite{Nezakati} considered the case where the covariance $\mbox{cov}(a_i, a_j) = 1- |i-j|/n$ and obtained that $E_n = O(\log n)^{1/2}$ for $n \to \infty$. 

Rezakhah and Shemehsavar \cite{Rezakhah1, Rezakhah2} studied random polynomials $Q_n(x) = \sum_{i=0}^{n-1} A_i x^i$ where the coefficients $A_i$ are successive Brownian images $A_i = W(t_i)$, $t_0< t_1< \ldots$ where $\{W(t), t\geq 0\}$ is the standard Brownian motion and analysed the asymptotic behavior of the real roots of the equations $Q_n(x) = K$ and $Q_n(x) = K x$ for $K$ constant (depending on $n$).  

In this paper, we  consider random polynomials  where the coefficients are correlated Gaussian random variables defined as successive {\it increments} of a fractional Brownian motion. These polynomials appear as natural generalisations of polynomials with independent coefficients. We study the asymptotic behaviour of the expected number of their real zeros and obtain unexpected features  different from previous studies which extend Kac's results in a remarkable way.   

 Recall that a fractional Brownian motion (FBM) is a family of Gaussian random variables $X =\{X(t): t\geq 0\}$ of mean 0  such that $X(0) = 0$, for any $0< t_1 < t_2 \ldots t_n$, the random vector $(X(t_1), X(t_2), \ldots, X(t_n))$ is also Gaussian and there exists $0 < H < 1$ (called the Hurst index) such that for all $t, s >0$,
              $$\mathbb{E}[(X(t) -  X(s))^2] = |t - s|^{2H}.$$ 
The discussion of the existence of FBM and various properties can be found in Kahane \cite{Kahane} and Nourdin \cite{Nourdin_Ivan}.         
We shall consider the increments $$a_0 = X(1), a_1 = X(2) - X(1), \ldots, a_n = X(n+1) - X(n), \ldots$$ of a fractional Brownian motion and study the number of real zeros of the random  polynomial
         $$P_n(x) = a_0 + a_1 x + \ldots+ a_{n-1} x^{n-1}.$$      
 All the random variables $a_i$ have the standard normal distribution. In the case $H = 1/2$, the process $X$ is the classical Brownian motion and it is well known that its increments are independent. Then, in that case, we retrieve the classical random polynomials where the coefficients $a_i$ and i.i.d random variables.    
 However for the case $H\ne 1/2$, the increments $a_0, a_1, \ldots$ are not independent since the quantity 
    \begin{eqnarray} \label{ew32wswaws2}
g(i, j, H):= \mathbb{E}[a_i a_j] = -|j-i|^{2H} + \frac{1}{2}|j-i+1|^{2H} + \frac{1}{2}|j-i-1|^{2H}
\end{eqnarray} 
does not always vanish for  $i \ne j$.     
For each  fixed $a_i$, 
           $$\lim_{n\to \infty} \mathbb{E}[a_i a_n] = 0.$$ This implies that the random variables $a_0, a_1, \ldots$ are asymptotically independent. 
We shall obtain the asymptotic estimate
     $$E_n \sim \frac{1}{\pi} (1 + 2\sqrt{H(1-H)}) \log n,\,\,\,\, n\to \infty.$$  
The main results of the paper are the following two theorems (the proofs are given in sections 3 and 4 respectively).
\begin{theorem} \label{th01}
Assume that $X$ is a fractional Brownian motion of Hurst index $H$, $a_0 = X(1)$ and $a_k = X(k+1) - X(k)$, $k = 1, 2, \ldots$ Then the average $E_n$ of the number of real zeros of the polynomial $P_n(x) = a_0 + a_1 x + a_2 x^2 + \ldots + a_{n-1} x^{n-1}$ is such that 
    $$E_n \sim \frac{1}{\pi} (1 + 2\sqrt{H(1-H)}) \log n,\,\,\,\, n\to \infty.$$ 
    Moreover, the average of the number of positive real zeros is $\sim (1/\pi) (2\sqrt{H(1-H)}) \log n$ and the average of the number of negative zeros is $\sim(1/\pi) \log n$.
\end{theorem}  

\begin{theorem} \label{th02}
 For any sequence of Gaussian random variables $X(1), X(2), \ldots$ satisfying 
\begin{eqnarray*} 
\mathbb{E}[X(k)^2] = 1, \,\, \mathbb{E}[X(k)] = 0\,\, \mbox{ and } \mathbb{E}[(X(k) - X(j))^2] = 1,\,\,\forall k\ne j,
\end{eqnarray*}
 let 
  $a_0 = X(1)$ and $a_k = X(k+1)-X(k)$, $k = 1,2,\ldots$ Then the average $E_n$ of the number of real zeros of the random polynomial $P_n(x) =a_0 + a_1 x + a_2 x^2 + \ldots + a_{n-1}x^{n-1}$ satisfies the asymptotic estimate
     $$E_n \sim (1/\pi) \log n,\,\,\, n\to \infty.$$
 Moreover the average of the number of positive zeros $E_n^+$ is such that $$\lim_{n \to \infty} E_n^+ = \frac{1}{3}  -\frac{ \log(2-\sqrt{3})}{\pi}.$$  
 \end{theorem}

The expression of $E_n$ reduces to $E_n \sim (2/\pi) \log n$ for $H = 1/2$, the formula obtained by Kac. It is surprising that the average number of negative zeros of $P_n(x)$ is independent of $H$ but the average number of positive zeros depends on $H$. 
Moreover the random polynomials corresponding to fractional Brownian motions of parameters $H$ and $1-H$ have on average, essentially, the same number of real zeros. This  was unexpected since these two processes have very different properties. In fact, this is the only non-obvious property depending on the Hurst index shared by these two processes that is known to the author. 

Theorem 2 extends the study to the limit case of $H = 0$. Although FBM of index $H = 0$ is not well defined, one can construct a sequence of standard Gaussian random variables $(a_i)$ such that (\ref{ew32wswaws2}) holds for $H = 0$. Then we shall prove that in that case, $E_n \sim (1/\pi) \log n$ and that the number of positive zeros is bounded in $n$. 
For the case of independent coefficients, the real zeros cluster symmetrically around $1$ and $-1$. However for $H\ne 1/2$, although the real zeros still cluster around $1$ and $-1$, they do so asymmetrically: the density at $-1$ is preserved while the density at $1$ decreases as $|H-1/2|$ increases so that it vanishes for the limit case $H=0$. The paper contains  some lengthy computations and all where checked using the Mathematica software. 
      
  \section{Kac - Rice formula}
Let $P_n(x) = a_0 + a_1 x + \ldots+ a_{n-1} x^{n-1}$ where $a_i = X(i+1) - X(i)$, $i = 0, 1, 2, \ldots$ and $X$ is a fractional Brownian motion of index $H$. 
Then by Kac-Rice formula (see for example Farahmand \cite{Farahmand} and  Bharucha-Ried and Sambandham \cite{Bharucha-Reid_Sambandham}), the expected number of real zeros of $P_n(x)$ is given by
 \begin{eqnarray}
 \label{kacgenf}
 E_n = \frac{1}{\pi} \int_{-\infty}^{\infty} \frac{\Delta^{1/2}}{\alpha}\ dx
 \end{eqnarray}
 where
  $$\alpha = \mathbb{E}[P_n(x)^2], \,\, \,\, \beta = \mathbb{E}[P_n'(x)^2],\,\,\, \gamma = \mathbb{E}[P_n(x) P_n'(x)] \,\, \mbox{ and } \Delta = \alpha \beta - \gamma^2.$$
As indicated by Kac himself and many others, the computation of this integral even in the particular case where $H = 1/2$ is very difficult. It is even more difficult in the presence of correlated random variables. The rest of the paper is devoted to finding an asymptotic estimate of  integral (\ref{kacgenf}) in terms of the Hurst index $H$ for $n\to \infty$. 

We shall explicitly compute the quantities $\alpha$, $\beta$, $\gamma$, $\Delta$ as functions of $H$, $n$ and $x$ as 
 \begin{eqnarray*} 
 \alpha & = & \sum_{i,j = 0}^{n-1} g(i, j, H)\ x^{i + j}\\
 \beta & =  & \sum_{i, j=1}^{n-1}  g(i, j, H) \ i \ j \ x^{i+j-2} \\
 \gamma & = & \sum_{i=0, j=1}^{n-1} g(i, j, H) \ j \ x^{i + j -1}.
 \end{eqnarray*}
We can make use of the immediate fact that $$ g(i, i, H) = 1, \,\, g(i, j, H) = g(j, i, H) = g(1, j-i+1, H) \mbox{ for }i\leq j$$ to obtain
  \begin{eqnarray} \label{eqsadw242}
 \alpha & = & \alpha_0 + 2 \sum_{j=2}^{n} g(1, j, H) \alpha_{j-1} \\ 
  \beta & =  & \beta_1 + 2 \sum_{j=2}^{n}  \ g(1, j, H) \beta_j \\
 \gamma & = & \gamma_0 + \sum_{j=2}^{n} g(1, j, H) \gamma_{j-1}
\end{eqnarray}
  \mbox{ with }\\
  \begin{eqnarray*}
  \alpha_j & =&  \sum_{k=0}^{n-j-1} x^{2k +j},\,\, \,\,j = 0, 1, 2, \ldots, n-1\\
\beta_j & = & \sum_{k=1}^{n-j} k(k+j-1) x^{2k+j-3},\,\,\,\, j= 1, 2, \ldots, n-1\\
\gamma_j & = & j x^{j-1} + \sum_{k = j+1}^{n-1} (2k -j) x^{2k-j-1},\,\, \,\, j=1,2,\ldots, n-1\\
\gamma_0 & = & \sum_{k =1}^{n-1} k \ x^{2k-1}\\
\beta_n & = &  0.
   \end{eqnarray*}

\section{Asymptotic estimate of $E_n$ for $0< H < 1$}
For $H = 1/2$, $\alpha = \alpha_0$, $\beta = \beta_1$ and $\gamma = \gamma_0$. We shall denote by $\Delta_0$ the function $\Delta$ corresponding to $H = 1/2$, that is, 
      $$\Delta_0 = \alpha_0 \beta_1 - \gamma_0^2.$$
Write  \begin{eqnarray*}
\Delta & = & \alpha \beta - \gamma^2 \\
& = & (\alpha_0 \beta_1  - \gamma_0^2) + 2 \sum_{j=2}^{n}  g(1, j, H) \left(\alpha_0  \beta_j + \beta_1 \alpha_{j-1} -\gamma_0 \gamma_{j-1}\right) \\
&&  + 4 \left(\sum_{j=2}^{n}  g(1, j, H) \beta_j\right) \left(\sum_{j=2}^{n}  g(1, j, H) \alpha_{j-1}\right) \\
&& - \left(\sum_{j=2}^{n}  g(1, j, H) \gamma_{j-1}\right)^2\\
& = & \Delta_0 + 2 \sum_{j=2}^{n}  g(1, j, H) \left(\alpha_0  \beta_j + \beta_1 \alpha_{j-1} -\gamma_0 \gamma_{j-1}\right) \\
&& + \sum_{j=2}^{n}  g(1, j, H)^2 (4\beta_j \alpha_{j-1} - \gamma_{j-1}^2) \\
&& + \sum_{2 \leq j < k \leq n} g(1, j, H) g(1, k, H) (4\beta_j \alpha_{k-1} + 4 \beta_k \alpha_{j-1} - 2 \gamma_{j-1}\gamma_{k-1}).
\end{eqnarray*}
We shall first obtain some asymptotic approximations, for $n$ large, of  each of the terms $\Delta_0$, \\$\alpha_0  \beta_j + \beta_1 \alpha_{j-1} -\gamma_0 \gamma_{j-1}$, $4\beta_j \alpha_{j-1} - \gamma_{j-1}^2$ and $(4\beta_j \alpha_{k-1} + 4 \beta_k \alpha_{j-1} - 2 \gamma_{j-1}\gamma_{k-1})$ by comparing them to the quantity $\alpha_0^2$. 
Let 
 \begin{eqnarray*} 
  A_j & = &(\alpha_0  \beta_j + \beta_1 \alpha_{j-1} -\gamma_0 \gamma_{j-1})/\alpha_0^2, \\
  B_j & = & (4\beta_j \alpha_{j-1} - \gamma_{j-1}^2)/\alpha_0^2,\\
  C_{j,k} & = & (4\beta_j \alpha_{k-1} + 4 \beta_k \alpha_{j-1} - 2 \gamma_{j-1}\gamma_{k-1})/\alpha_0^2.
  \end{eqnarray*}
Then 
  \begin{eqnarray*}
   \frac{\Delta}{\alpha^2} = \frac{(\Delta_0/\alpha_0^2) + 2\sum_{j=2}^n g(1, j, H) A_j + \sum_{j=2}^n g(1,j,H)^2 B_j + \sum_{2 \leq j < k \leq n} g(1, j, H) g(1, k, H) C_{j,k}}{\left (1 + 2 \sum_{j=2}^n g(1, j, H) (\alpha_{j-1}/\alpha_0)\right)^2}.
  \end{eqnarray*}
 Direct computations yield the following:\\ 
    \begin{eqnarray*}
 \frac{\Delta_0}{\alpha_0^2}  & = &  \frac{(1 + 2 (-1 + n^2) x^{2 n} + x^{4 n} - n^2 x^{2 n-2} - 
 n^2 x^{2 n+2})}{(x^2-1)^2 (x^{2 n}-1)^2},\\ 
  A_j& = & \frac{x^{-3 - j} (2 x^{2 + 2 j} + 2 x^{4 + 4 n} - 
   x^{2 n} (x^2 + x^{2 j}) (2 x^2 + 
      n (1 - j + n) (x^2-1)^2))}{( x^2-1)^2 (x^{2 n}-1)^2}, \\
B_j &  = &  \frac{1}{(x^2-1)^2 (x^{
   2 n} -1)^2}   \left[x^{-2 j} (-((j-3) x^{
       2 j} - (j-1) x^{-2 + 2 j} + (1 - j + 2 n) x^{2 n} \right. \\
&&  + \left. (1 + j - 2 n) x^{2 + 2 n})^2 
 - 4 (x^{2 j-2} - x^{2 n}) (-2 x^{2 + 2 j} + n (1 + n) x^{2 n}\right. \\
&& \left.  - 
      2 (n^2-1) x^{2 + 2 n} + (n-1) n x^{4 + 2 n} + 
      j (x^2-1) (x^{2 j} + 
         n x^{2 n} - (n-1) x^{2 + 2 n})))\right],
\end{eqnarray*}
and
 \begin{eqnarray*}
 C_{j,k} & = & \frac{2 x^{-j-k-4}}{(x^2-1)^2 (x^{2n}-1)^2} \left[x^{2 (j + k)} (-1 + k - 2 (-3 + k) x^2 + (-1 + k) x^4\right. \\
&& \left. -  j (-1 + k) (-1 + x^2)^2) + 
 x^{4 + 4 n} (-1 + k - 2 (-3 + k) x^2 + (-1 + k) x^4 \right. \\
 && \left. - j (-1 + k) (-1 + x^2)^2) - 
 x^{2 + 2 n} (x^{2 j} + x^{2 k}) (1 + 2 n (2 + n) + 2 x^2 - 
    4 n (2 + n) x^2 \right.\\
 && \left. + (1 + 2 n (2 + n)) x^4 + 
    j (-1 + k - 2 n) (-1 + x^2)^2 - k (1 + 2 n) (-1 + x^2)^2)\right].
  \end{eqnarray*}
Moreover,
$$\alpha_j/\alpha_0 = \frac{x^{-j} (x^{2 j} - x^{2 n})}{(1 - x^{2 n})}.$$
The integral of $\Delta^{1/2}/\alpha$ will be decomposed into integrals over the subsets (1) $|x|> 1+1/n$, (2) $1 < |x| < 1+1/n$ and (3) $|x| <1$. 
\paragraph  {\bf 1) Case $|x| > 1+1/n$.}  
  For $|x| > 1+1/n$, elementary computations give the following asymptotic approximations (for large $n$): 
          \begin{eqnarray*}
          \Delta_0/\alpha_0^2 & = & \frac{1}{(x^2-1)^2} - \frac{n^2 x^{2n -2}}{(x^{2n}-1)^2} =  \frac{1}{(x^2 - 1)^2} + x^{-2n}O(n^2), \\
          A_j & = & \frac{2 x^{1-j}}{(x^2-1)^2}  + \frac{2 x^{-n}}{(x^2-1)^2} + x^{-n} O(n^2),\\
          B_j & = & \frac{x^{-2j} W_1(x)}{(x^2-1)^2} + \frac{x^{-2n-4} W_1(x)}{(x^2-1)^2}  - \frac{2 x^{-2n-2} (x^2+1)^2}{(x^2-1)^2} + x^{-2n-2} O(n^2),  \\
            C_{j,k} & = & \frac{2 x^{-j-k}}{(x^2-1)^2} W_2(x) +\frac{2 x^{j - k - 2n -2}}{(x^2-1)^2} W_3(x) +  \frac{2 x^{k -j- 4 n -2}}{(x^2-1)^2} W_3(x) + \frac{2 x^{j+k - 4n -4}}{(x^2-1)^2} W_2(x) \\
            && \quad \quad + \frac{1}{(x^2-1)^2} O(x^{-4n}),
          \end{eqnarray*}
with 
 \begin{eqnarray*}
 W_1(x) & =  & -1 + 6 x^2 - x^4 + (2j-j^2) (x^2-1)^2,\\
 W_2(x)&  = & -1 + 6x^2 - x^4 + (j+k - jk) (x^2-1)^2,\\
 W_3(x) & = & -(j (-1 + k - 2 n) + 2 n (2 + n) - k (1 + 2 n)) (-1 + 
      x^2)^2 - (1 + x^2)^2).\\
\end{eqnarray*}     
Since $|x|>1+1/n$ implies that $(x^2 - 1)^2 > 4/n^2$, then (using $j, k \leq n$),
  \begin{eqnarray*}
          \Delta_0/\alpha_0^2 & = & \frac{1}{(x^2 - 1)^2} + x^{-2n}O(n^2),\\
          A_j & = & \frac{2 x^{1-j}}{(x^2-1)^2}  + x^{-n} O(n^2),\\
          B_j & = & \frac{x^{-2j} W_1(x)}{(x^2-1)^2} + x^{-2n+2} O(n^2),\\
            C_{j,k} & = & \frac{2 x^{-j-k}}{(x^2-1)^2} W_2(x) +x^{-2n} O(n^2).
                   \end{eqnarray*}
Moreover, 
        $$\alpha_j/\alpha_0 = x^{-j} + O(x^{-2n}).$$
Using 
 $$\frac{a}{b} = 1 + \frac{a-b}{b}$$ and the obvious fact that $|g(1, j, H)| \leq 1$, we obtain (after some elementary computations) that
\begin{eqnarray} \label{ewdwdsswwqewe}
\frac{\Delta}{\alpha^2} & =& \frac{1}{(x^2-1)^2}\left[1 - \frac{(x^2-1)^2 \left(\sum_{j=2}^{n} g(1, j, H) x^{-j} (j-1)\right)^2}{\left(1 + 2 \sum_{j=2}^{n} g(1, j, H) x^{1-j}\right)^2} \right] + x^{-n}O(n^4).
\end{eqnarray}
We shall rewrite the involved sums by using the classical  Lerch transcendent function
$$\Phi(z, s, a) = \sum_{k=0}^\infty z^k/(k+a)^s,\quad \quad \,\, |z| < 1.$$
From the definition of $g(1,j,H)$ (relation (\ref{ew32wswaws2})), it easy to obtain that
\begin{eqnarray*}
\sum_{j=2}^{n} g(1, j, H) x^{-j} (j-1) & = &(1/2)(1/x)^{n+1} (-\Phi(1/x, -1 - 2 H, n-1) + 2 
\Phi(1/x, -1 - 2 H, n))\\
&& + (1/2)(1/x)^{n+1} (-\Phi(1/x, -1 - 2 H,  n + 1) -   \Phi(1/x, -2 H, n-1))\\
  && + (1/2)(1/x)^{n+1} \Phi(1/x, -2 H, n+1) +   (1/2) \Phi(1/x, -1 - 2 H, 0) \\
&& + (1/2) x^{-2}\Phi(1/x, -1 - 2 H, 0)  - (1/x)\Phi(1/x,-1 - 2 H, 0)\\
&& - (1/2) \Phi(1/x,-2 H,0) +   (1/2)x^{-2}\Phi(1/x, -2 H,0),
 \end{eqnarray*}
and
  \begin{eqnarray*}
\sum_{j=2}^{n} g(1, j, H) x^{1-j}& = &(1/2)x^{-1} \Phi(1/x, -2 H, 2) - (1/2)x^{-n}\Phi(1/x, -2 H, n-1)\\
&& + x^{-n}\Phi(1/x, -2 H, n) - (1/2)x^{-n}\Phi(1/x, -2 H, n+1) \\
&& + (1/2)x^{-1}\Phi(1/x, -2H, 0)  - \Phi(1/x, -2H, 0).
\end{eqnarray*} 
For all $|z|< 1$, $a>0$, $s<0$, 
   $$|\phi(z, s, a)| \leq \sum_{k=0}^\infty \frac{|z|^k}{a^s} = \frac{a^{-s}}{1-|z|}.$$ Then for $s$ fixed $$\Phi(1/x, s, m) = \frac{1}{1-|1/x|}O(m^{-s}),\,\,\,m \to \infty.$$
For example, for all $|x| > 1+1/n$, 
    $$|\Phi(1/x, -1 - 2 H, n-1)| \leq \frac{(n-1)^{2H+1}}{1-|1/x|} \leq (n+1)(n-1)^{2H+1} $$ and hence
         $$\Phi(1/x, -1 - 2 H, n-1)  = O(n^{2H+2}).$$ 
This yields
        \begin{eqnarray*}
\sum_{j=2}^{n} g(1, j, H) x^{-j} (j-1) & = & (2x^2)^{-1}(x-1)\left[(x-1)\Phi(1/x, -1-2H,0) -(1+x)\Phi(1/x, -2H,0)\right]\\
&& + x^{-n-1}O(n^{2H+2}),\\
\sum_{j=2}^{n} g(1, j, H) x^{1-j}& = & -\frac{1}{2} + \frac{(x-1)^2}{2x}\Phi(1/x, -2H,0) +  x^{-n} O(n^{2H+1})
\end{eqnarray*}
where we have used the identity
 $$\Phi(1/x, -2 H, 2) = x^2\Phi(1/x, -2 H, 0)- x.$$
  Using the classical polylogarithm function
  $$\mbox{Li}_s(z) = \Phi(z, s, 0)= \sum_{k=1}^\infty z^k/k^s,\quad \quad \,\, |z| < 1, $$ yields
       \begin{eqnarray*}
  \sum_{j=2}^n g(1, j, H) x^{-j} (j-1) & = & (2x^2)^{-1}(x-1)\left[(x-1)\mbox{Li}_{s}(1/x) -(1+x)\mbox{Li}_{t}(1/x))\right] + x^{-n-1}O(n^{2H+2}), \\
  \sum_{j=2}^n g(1, j, H) x^{1-j} & = & -\frac{1}{2} + \frac{(x-1)^2}{2x}\mbox{Li}_t (1/x) + x^{-n} O(n^{2H+1})
  \end{eqnarray*}
with $s = -1-2H$ and $t = -2H$.\\
Then (\ref{ewdwdsswwqewe}) implies
     \begin{eqnarray*}
\frac{\Delta}{\alpha^2} & =  &\frac{1}{(x^2-1)^2}\left(1 - \frac{(1+x)^2\left[(-1+x)\mbox{Li}_{s}(1/x) -(1 + x) \mbox{Li}_{t}(1/x)\right]^2}{4 x^2 \mbox{Li}_t(1/x)^2}\right)\\
&& \qquad + x^{-n} O(n^{2H+2}) +  x^{-n} O(n^{4}).
\end{eqnarray*} 
Therefore 
   \begin{eqnarray*}
\frac{\Delta^{1/2}}{\alpha} &= & \frac{1}{|x^2-1|}\left(1 - \frac{(1+x)^2\left[(-1+x)\mbox{Li}_{s}(1/x) -(1 + x) \mbox{Li}_{t}(1/x)\right]^2}{4 x^2 \mbox{Li}_t(1/x)^2}\right)^{1/2}\\
&& \qquad + x^{-n} O(n^{2H+2}) +  x^{-n} O(n^{4}).
\end{eqnarray*} 
For any fixed $0<\epsilon < 1$, we have that the sequences
$$n \to \int_{1+ n^{\epsilon - 1}} ^\infty n^{4} x^{-n} dx \,\,\, \mbox{ and } \,\,\, n \to \int_{-\infty}^{-1- n^{\epsilon - 1}} n^{4} x^{-n} dx$$ are bounded in $n$ (in fact both sequences converge to 0). 
Then for any  $0<\epsilon < 1$ fixed, 
      \begin{eqnarray*}
\int_{|x| \geq  1+ n^{\epsilon - 1}} \frac{\Delta^{1/2}}{\alpha} dx &= & \int_{|x|\geq 1+ n^{\epsilon - 1}} \frac{1}{|x^2-1|}\ell(x)^{1/2} dx + O(1)
\end{eqnarray*}
where $\ell$ is the function defined for $|x|>1$ by
   \begin{eqnarray} \label{wewsdw23ell}
\ell(x) = \left(1 - \frac{(1+x)^2\left[(-1+x)\mbox{Li}_{s}(1/x) -(1 + x) \mbox{Li}_{t}(1/x)\right]^2}{4 x^2 \mbox{Li}_t(1/x)^2}\right).
\end{eqnarray}
 The function $\ell(x)$ is continuous and bounded in $(-\infty, -1) \cup (1, \infty)$. Moreover, it is increasing on each of these intervals (separately) and admits the following limits:
  \begin{eqnarray*}
  \lim_{x \searrow  1} \ell(x)^{1/2} & = & C(H) = (4H(1-H))^{1/2}\\
  \lim_{x \nearrow  -1} \ell(x)^{1/2} & = &  1 \\
  \lim_{x\to \pm \infty} \ell(x)^{1/2} & =& M(H) = (2^{-2 (1 + s + t)} (-2^s + 2^t) (2^s - 2^t + 2^{2 + s + t}))^{1/2}
  \end{eqnarray*}
where $s = -1-2H$ and $t = -2H$. \\
Then \begin{eqnarray*}
\int_{1+ n^{\epsilon - 1}}^{\infty}  \frac{1}{x^2-1} \ell(x)^{1/2}dx & \leq & M(H) \int_{1+ n^{\epsilon - 1}}^{\infty} \frac{1}{x^2-1} dx\\
&=  & M(H) (i \pi/2 + \arctanh(1 + n^{\epsilon -1})\\
 & =  & M(H)\frac{(1-\epsilon)}{2} \log n + O(1).
 \end{eqnarray*}  
Also, 
     \begin{eqnarray*}
\int_{1+ n^{\epsilon - 1}}^{\infty}  \frac{1}{x^2-1} \ell(x)^{1/2} dx& \geq & C(H) \int_{1+ n^{\epsilon - 1}}^{\infty} \frac{1}{x^2-1} dx\\
 & =  & C(H)\frac{(1-\epsilon)}{2} \log n + O(1).
 \end{eqnarray*}  
We want to show that 
   \begin{eqnarray*}
\int_{1+ n^{\epsilon - 1}}^{\infty}  \frac{1}{x^2-1} \ell(x)^{1/2}dx = C(H)\frac{(1-\epsilon)}{2} \log n + O(1).
\end{eqnarray*}
For any small fixed number $\delta>0$ and $\eta = \ell^{-1}(\delta + C(H))$, it is the case that 
 $$\int_{1+ n^{\epsilon - 1}}^{\eta}  \frac{1}{x^2-1} \ell(x)^{1/2}dx \leq (C(H) + \delta)\int_{1+ n^{\epsilon - 1}}^{\eta}  \frac{1}{x^2-1} dx$$
because the function $\ell(x)$ is increasing. Moreover since for any fixed $\eta >1$, 
      $$\int_\eta^\infty  \frac{1}{x^2-1} \ell(x)^{1/2}dx < \infty$$ it follows that
\begin{eqnarray*}
 \int_{1+ n^{\epsilon - 1}}^{\infty}  \frac{1}{x^2-1} \ell(x)^{1/2}dx &  = & \int_{1+ n^{\epsilon - 1}}^{\eta} \frac{1}{x^2-1} \ell(x)^{1/2}dx + \int_\eta^\infty  \frac{1}{x^2-1} \ell(x)^{1/2}dx \\
 & \leq  & (C(H) + \delta)\int_{1+ n^{\epsilon - 1}}^{\eta}  \frac{1}{x^2-1} dx + O(1). 
\end{eqnarray*}
Since $\delta$ can be taken arbitrarily small, this implies that, for all large $n$, 
     \begin{eqnarray*}
 \int_{1+ n^{\epsilon - 1}}^{\infty}  \frac{1}{x^2-1} \ell(x)^{1/2}dx = C(H)\int_{1+ n^{\epsilon - 1}}^{\eta}  \frac{1}{x^2-1} dx + O(1) = C(H)\frac{(1-\epsilon)}{2} \log n + O(1).  
\end{eqnarray*}
Because the number $\epsilon$ can be chosen arbitrary small, then
   $$\int_{1+ n^{-1}}^{\infty}  \frac{1}{x^2-1} \ell(x)^{1/2}dx  = C(H)\frac{1}{2} \log n + O(1)  =  (H(1-H))^{1/2} \log n + O(1).$$ 
     
This can be repeated for the integral on $(-\infty, -1-1/n]$ by using the fact that on this interval, $\ell(x)$ is increasing and bounded by 1. Then
\begin{eqnarray*}
\int_{-\infty}^{-1- n^{\epsilon - 1}}  \frac{1}{x^2-1} \ell(x)^{1/2}dx & \leq &  \int_{-\infty}^{-1- n^{\epsilon - 1}} \frac{1}{x^2-1} dx\\
&=  &\frac{(1-\epsilon)}{2} \log n + O(1).
\end{eqnarray*}
We can also only consider values around $x = -1$ as previously and obtain that
   \begin{eqnarray*}
\int_{-\infty}^{-1- n^{\epsilon - 1}} \frac{1}{x^2-1} \ell(x)^{1/2}dx  = \frac{(1-\epsilon)}{2} \log n + O(1).
\end{eqnarray*}
Because $\epsilon$ can be taken arbitrarily small, this implies 
           \begin{eqnarray*}
\int_{-\infty}^{-1- n^{\epsilon - 1}} \frac{1}{x^2-1} \ell(x)^{1/2}dx  = (1/2) \log n + O(1).
\end{eqnarray*}
Then we obtain that 
\begin{eqnarray} \label{EQ01}
  \int_{|x| > 1+ n^{-1}} \left(\frac{\Delta^{1/2}}{\alpha}\right) dx = (1/2 + \sqrt{H(1-H)}) \log n + O(1). 
  \end{eqnarray}
        
\paragraph  {\bf 2) Case $1< x < 1+1/n.$} 
For the values of $x$ near $\pm 1$,  we make use of Taylor's expansion.\\ 
 a) First,  
\begin{eqnarray*}
\frac{\Delta_0}{\alpha_0^2} & =& \frac{x^2 - n^2 x^{2 n} + 2 (-1 + n^2) x^{2 + 2 n} - n^2 x^{4 + 2 n} + x^{
 2 + 4 n}}{x^2 (x^2-1)^2(x^{2 n} - 1)^2}.
 \end{eqnarray*}
 Let $$h_n = (x^2-1)^2 \Delta_0/\alpha_0^2 = \frac{x^2 - n^2 x^{2 n} + 2 (-1 + n^2) x^{2 + 2 n} - n^2 x^{4 + 2 n} + x^{
 2 + 4 n}}{x^2 (x^{2 n} - 1)^2}.$$
Write $$p_n(x) = x^2 - n^2 x^{2 n} + 2 (-1 + n^2) x^{2 + 2 n} - n^2 x^{4 + 2 n} + x^{
 2 + 4 n}$$  
(the numerator of $h_n(x)$). Then elementary computations show that
 $$p_n(\pm 1) = p_n'(\pm 1) = p_n''(\pm 1) = p_n'''(\pm 1) = 0 \mbox{ and } p_n^{(4)}(x) = O(n^4), \mbox{ for all } 1 < |x|< 1+1/n.$$
In fact
 for $x = \pm (1 + c/n)$ with $0 \leq c \leq 1 $,
 $$ \lim_{n\to \infty} \frac{p_n^{(4)}(x)}{n^4} =(32 e^{2/c} (-2 - 8 c + c^2 (-7 + 8 e^{2/c})))/c^2.$$
 Then 
    $$p_n(x) = (x-1)^4 O(n^4) \mbox{ and } p_n(x) = (x+1) O(n^4)  \mbox{ for } n \mbox{ large}. $$
Moreover noting that, for $1 < |x| < 1 + 1/n$,
                $$(x^{2n} - 1)^2 = (|x|^{2n}-1)^2 =  (|x|-1)^2(1 + |x| + |x|^2 + \ldots + |x|^{2n-1})^2 >  (2n)^2(|x|-1)^2, $$ 
                which implies that 
    $$(x^{2n}-1)^2 > 4n^2 (x-1)^2 \mbox{ for } 1 < x < 1+1/n $$
 and $$ (x^{2n}-1)^2 > 4n^2 (x+1)^2 \mbox{ for } -1 > x > -(1+1/n),$$ then            
   \begin{eqnarray} \label{wewsqa1}
h_n(x)&= &  \frac{(x-1)^4 O(n^4)}{x^2 (x-1)^2 n^2} =  (x-1)^2 O(n^2)\,\, \mbox{ for } 1 < x < 1+1/n,\\
h_n(x)&= & \frac{(x+1)^4 O(n^4)}{x^2 (x+1)^2 n^2} =  (x+1)^2 O(n^2)\,\, \mbox{ for } -1 > x > -(1+1/n).
\end{eqnarray}     
It follows that 
  \begin{eqnarray*}
\frac{\Delta_0}{\alpha_0^2} & = &  \frac{1}{(x+1)^2} O(n^2),\,\,\,1 < x < 1+1/n\\
 \frac{\Delta_0}{\alpha_0^2} & = & \frac{1}{(x-1)^2} O(n^2),\,\, -(1+1/n) < x < -1.
\end{eqnarray*}  
(b) Write
    \begin{eqnarray*}
A_j=\frac{\alpha_0  \beta_j + \beta_1 \alpha_{j-1} -\gamma_0 \gamma_{j-1}}{\alpha_0^2} = \frac{x^{-j}(x^2-1)^2 t_n(x)}{(x^{2n} - 1)}.
\end{eqnarray*}
Then elementary computations yield
       $$t_n(\pm 1) = t_n'(\pm 1) = t_n''(\pm 1) = t_n'''(\pm 1) = 0 \mbox{ and } t_n^{(4)}(x) = O(n^4), \mbox{ for all } 1 < |x| < 1+1/n.$$ 
Then 
    $$t_n(x) = (x-1)^4 O(n^4) \mbox{ and } t_n(x) = (x+1)^4 O(n^4) \mbox{ for } n \mbox{ large} $$  and as previously, 
\begin{eqnarray*}
A_j & = & \frac{x^{-j}}{(x+1)^2} O(n^2), \,\,1 < x < 1+1/n,\\
A_j & = & \frac{x^{-j}}{(x-1)^2} O(n^2),\,\,-(1+ 1/n) < x < -1.
\end{eqnarray*}
(c) Similar computations yield that 
        \begin{eqnarray*}
 &&B_j = \frac{4\beta_j \alpha_{j-1} - \gamma_{j-1}^2}{\alpha_0^2} = \frac{x^{-2j}}{(x+1)^2} O(n^2),\,\, 1 < x < 1+1/n,\\
  &&B_j=\frac{4\beta_j \alpha_{j-1} - \gamma_{j-1}^2}{\alpha_0^2}= \frac{x^{-2j}}{(x+1)^2} O(n^2),\,\, -(1+1/n) < x < -1,\\
  \end{eqnarray*}
 (d)
  \begin{eqnarray*}
 &&C_{j,k} = \frac{4\beta_j \alpha_{k-1} + 4 \beta_k \alpha_{j-1} - 2 \gamma_{j-1}\gamma_{k-1}}{\alpha_0^2}
 =\frac{x^{-j-k}}{(x+1)^2} O(n^2),\,\, \,1 < x < 1+1/n,\\
 &&C_{j,k} = \frac{4\beta_j \alpha_{k-1} + 4 \beta_k \alpha_{j-1} - 2 \gamma_{j-1}\gamma_{k-1}}{\alpha_0^2}
 =\frac{x^{-j-k}}{(x+1)^2} O(n^2), \,\,\, -1-1/n < x < -1.
 \end{eqnarray*}
 Then, for $1 < x < 1+1/n$,
     \begin{eqnarray*}
\frac{\Delta}{\alpha^2} & =& \left(\frac{\Delta}{\alpha_0^2} \right) \left(\frac{\alpha_0^2}{\alpha^2}\right) \\
& = & \left(\frac{\Delta}{\alpha_0^2}\right) \left(\frac{1}{1 + 2\sum_{j=2}^n g(1, j, H) (\alpha_{j-1}/\alpha_0)}\right)^2,
\end{eqnarray*}
where 
\begin{eqnarray*}
\Delta/\alpha_0^2 & = & \frac{1}{(x+1)^2} O(n^2) + 2 \sum_{j=2}^{n}  g(1, j, H) \frac{x^{-j}}{(x+1)^2} O(n^2)\\
&& + \sum_{j=2}^{n}  g(1, j, H)^2 \frac{x^{-2j}}{(x+1)^2} O(n^2) \\
&& + \sum_{2 \leq j < k \leq n} g(1, j, H) g(1, k, H) \frac{x^{-j-k}}{(x+1)^2} O(n^2)
\end{eqnarray*}
and as previously $$\alpha_j/\alpha_0 = x^{-j} + O(x^{-2n}).$$
 Then  for all $1 < x < 1+1/n$,        
  \begin{eqnarray*}
\frac{\Delta}{\alpha^2} & = & \frac{1}{(x+1)^2} \frac{(1 + \sum_{j=2}^{n}  g(1, j, H) x^{-j})^2}{(1 + \sum_{j=2}^{n}  g(1, j, H) x^{1-j})^2} O(n^2) \\
&= & \frac{1}{(x+1)^2} O(n^2).
\end{eqnarray*}
Hence
 \begin{eqnarray*}
\frac{\Delta^{1/2}}{\alpha} & = & \frac{1}{x+1} O(n) 
\end{eqnarray*}
 and therefore,
              \begin{eqnarray} \label{EQ02}
 \int_{1}^{1+1/n} \frac{\Delta^{1/2}}{\alpha} dx = O(1).
\end{eqnarray}      
Similarly, for $-(1+1/n) < x < -1$,
     \begin{eqnarray*}
\frac{\Delta}{\alpha^2} & = & \frac{1}{(x-1)^2} O(n^2) 
\end{eqnarray*}
 and hence  \begin{eqnarray} \label{EQ03}
 \int_{-(1+1/n)}^{-1} \frac{\Delta^{1/2}}{\alpha} dx = O(n) \int_{-(1+1/n)}^{-1} \frac{1}{|x-1|}dx = O(1).
 \end{eqnarray}
 We conclude from (\ref{EQ01}), (\ref{EQ02}) and (\ref{EQ03}) that
 $$   \int_{|x| > 1} \left(\frac{\Delta^{1/2}}{\alpha}\right) dx = (1/2 + \sqrt{H(1-H)}) \log n + O(1) $$
with  
     \begin{eqnarray*}
     \int_1^\infty \frac{\Delta^{1/2}}{\alpha} dx & = & \sqrt{H(1-H)}) \log n + O(1),  \\
      \int_{-\infty}^{-1} \frac{\Delta^{1/2}}{\alpha} dx & = & (1/2) \log n + O(1).
      \end{eqnarray*}    

\paragraph  {\bf 3)  Case $|x| < 1$.}
For any fixed $|x|< 1$, (with $x\ne 0$), it is easy to obtain the following estimates for large $n$: \\
 \begin{eqnarray*}
\frac{\alpha_j}{\alpha_0} & = & x^{j} + O(x^{2n-j}),\\
  \frac{\beta_j}{\beta_1} & = & \frac{x^{j-1}(j + (2-j)x^2)}{x^2 + 1} + O(n^2 x^{2n -j}),\\
  \frac{\gamma_j}{\gamma_0} & = &  x^{j-2} (j+ (2-j)x^2)  + O(n x^{2n -j}).
\end{eqnarray*}      
 Now write for $|x|< 1$, $j=1,2,\ldots, n-1,$
    $$\phi_j(x) =   x^{j},\,\,\, \xi_j(x) = \frac{x^{j-1}(j + (2-j)x^2)}{1 + x^2},  \,\,\eta_j(x) = (j + 2x^2 -j x^2) x^{j-2}\mbox{ and } \xi_n(x) = 0.$$
Then, for large $n$ (using the fact that $j\leq n$ and $|g(1,j,H)| \leq 1$)  
   \begin{eqnarray*}
   \frac{\alpha}{\alpha_0} & = & 1 + 2\sum_{j=2}^{n} g(1, j, H) \phi_{j-1}(x) + O(n x^{n})\\
   \frac{\beta}{\beta_1}&  = &  1 + 2 \sum_{j=2}^{n}  g(1, j, H) \xi_j(x) + O(n^3 x^{n})\\
   \frac{\gamma}{\gamma_0}  & = & 1 + \sum_{j=2}^{n}  g(1, j, H)  \eta_{j-1}(x) + O(n^2 x^{n}). 
   \end{eqnarray*} 
Therefore, 
\begin{eqnarray*}
\frac{\Delta}{\alpha^2} & = &  \frac{\beta}{\alpha} - \frac{\gamma^2}{\alpha^2}\\
  & = & \frac{\beta_1}{\alpha_0} \left(\frac{1 + 2 \sum_{j=2}^{n}  g(1, j, H) \xi_j(x)}{1 + 2\sum_{j=2}^{n} g(1, j, H) \phi_{j-1}(x)}\right) - \frac{\gamma_0^2}{\alpha_0^2} \left(\frac{1 + \sum_{j=2}^{n}  g(1, j, H)  \eta_{j-1}(x)}{1 + 2\sum_{j=2}^{n} g(1, j, H) \phi_{j-1}(x)}\right)^2 \\
  &&\quad \qquad +\left(\frac{\beta_1}{\alpha_0}- \frac{\gamma_0^2}{\alpha_0^2}\right) O(n^3\ x^{n}).
 \end{eqnarray*}
For $|x| <1$, unlike the case $|x|>1$, the quantities $\beta_1/\alpha_0$ and $\gamma_0/\alpha_0$ are bounded in $n$ and in fact, 
 $$\frac{\beta_1}{\alpha_0} = \frac{x^2+1}{(x^2-1)^2} + O(n^2 x^{2n})\,\,\, \mbox{ and } \left(\frac{\gamma_0}{\alpha_0}\right)^2 = \frac{x^2}{(x^2 - 1)^2} + O(n x^{2n}).$$  
So we write
 \begin{eqnarray*}
\frac{\Delta}{\alpha^2} & = & \frac{x^2+1}{(x^2-1)^2} \left(\frac{1 + 2 \sum_{j=2}^{n}  g(1, j, H) \xi_j(x)}{1 + 2\sum_{j=2}^{n} g(1, j, H) \phi_{j-1}(x)}\right)\\
&&\quad \qquad \qquad  - \frac{x^2}{(x^2 - 1)^2} \left(\frac{1 + \sum_{j=2}^{n}  g(1, j, H)  \eta_{j-1}(x)}{1 + 2\sum_{j=2}^{n} g(1, j, H) \phi_{j-1}(x)}\right)^2 + O(n^5 \ x^n).
 \end{eqnarray*}
 After some elementary transformations, this yields
  \begin{eqnarray*}
\frac{\Delta}{\alpha^2} & = &  \frac{x^2+1}{(x^2-1)^2} \left(\frac{1 + A_1}{1 + B_1}\right)
- \frac{x^2}{(x^2 - 1)^2} \left(\frac{1 + C_1}{1 + B_1}\right)^2 + O(n^5 \ x^n)
\end{eqnarray*}
with 
  $$A_1 = 2 \sum_{j=2}^{\infty}  g(1, j, H) \xi_j(x),\,\, B_1 = 2\sum_{j=1}^{\infty} g(1, j, H) \phi_{j-1}(x),\,\, C_1 = \sum_{j=2}^{\infty}  g(1, j, H)  \eta_{j-1}(x).$$
Now, as previously, these quantities can be represented in terms of the polylogarithm function as follows:
 \begin{eqnarray*}
 A_1 & = & -\frac{((-1 + x)^3) (1 + x) \mbox{Li}_s(x) + 
  x (1 + x^2 + 2 (-1 + x)^2 \mbox{Li}_t(x)))}{(x + x^3)} \\
 B_1 & = & -1 + (-2 + x^{-1} + x) \mbox{Li}_t(x) \\
 C_1 & = & \frac{-2 x^3 - (-1 + x)^3 (1 + x) \mbox{Li}_s(x)+ (-1 + x)^2 (-1 - 2 x + x^2) \mbox{Li}_t(x)}{2 x^3}
 \end{eqnarray*}
 where, as before, $s = -1-2H$ and $t = -2H$. 
Elementary computations show that
      $$(x^2+1) \left(\frac{1 + A_1}{1 + B_1}\right)
- x^2\left(\frac{1 + C_1}{1 + B_1}\right)^2 = \ell(1/x)$$ where 
 $\ell$ is the function defined by (\ref{wewsdw23ell}). 
 Hence
    $$\frac{\Delta}{\alpha^2} = \frac{\ell(1/x)}{(x^2-1)^2} + O(n^5 x^n),\,\,\,\, -1 < x < 1.$$ 
We can now compute the integral of $\Delta^{1/2}/\alpha$ on the intervals $[0, 1), (-1, 0]$ as previously by first integrating over $[0, 1-1/n)$ and $(-1+1/n, 0]$ and obtain that 
  \begin{eqnarray*}     
   \int_0^1 \left(\frac{\Delta^{1/2}}{\alpha}\right) dx = \sqrt{H(1-H)} \log n + O(1)\\
   \int_{-1}^0 \left(\frac{\Delta^{1/2}}{\alpha}\right) dx  = (1/2) \log n + O(1). 
  \end{eqnarray*}
We conclude that 
      $$\int_{-\infty}^{\infty} \left(\frac{\Delta^{1/2}}{\alpha}\right) dx = (1 + 2\sqrt{H(1-H)})\log n + O(1)$$ and therefore
       $$E_n = \frac{1}{\pi} (1 + 2\sqrt{H(1-H)}) \log n + O(1).$$
 This concludes the proof of Theorem 1.         

\section{The limit cases $H=0$ and $H=1$} 
For $H = 1$, fractional Brownian motion reduces to a random linear function $X(t) = X(1) t$ (where $X(1)$ has the standard normal distribution). The corresponding random polynomial reduces to $$P_n(x) = X(1) (1 + x + x^2 + \ldots + x^{n-1})$$ and has at most one real root. This has no practical interest. Kac's formula is not applicable because the quantity $\Delta = 0$ for all values of $x$. 
When $H$  approaches 0, the fractional Brownian motion is not well defined. 
However it is easy to construct a sequence of Gaussian random variables $X(1), X(2), X(3),\ldots$ such that 
 \begin{eqnarray} \label{weww23rea12}
\mathbb{E}[X(k)^2] = 1, \,\, \mathbb{E}[X(k)] = 0\,\, \mbox{ and } \mathbb{E}[(X(k) - X(j))^2] = 1,\,\,\forall k\ne j
\end{eqnarray} and hence $$\mathbb{E}[X(k) X(j)] = 1/2 = (k^{2H}  + j^{2H} - |k-j|^{2H})/2\,\,\,\mbox{ with } H = 0.$$
To obtain such a sequence, consider a sequence $\{Z_i\}$ of independent standard Gaussian variables and iteratively define
$X(1) = Z_1$, $X(2) = c_1 Z_1 + c_2 Z_2$ where $c_1$ and $c_2$ are real numbers such that $\mathbb{E}[X(2) X(1)] = 1/2$ and $\mathbb{E}[(X(2)^2] = 1$, for example, $a = 1/2$ and $b = \sqrt{3}/2$. In general write
$X(n) = c_1 Z_1 + c_2 Z_2 + \ldots + c_n Z_n$ and compute the coefficients $c_k$ by using the relations
 $\mathbb{E}[X(n) X(k)] = 1/2$, $k = 1, 2, \ldots, n-1$ and $\mathbb{E}[(X(n)^2] = 1$. 
As previously, we consider the increments
    $$a_k = X(k+1) - X(k),\,\,k = 0, 1, 2, \ldots \mbox{ where }\,\, X(0) = 0.$$ Clearly, $\mathbb{E}[a_k^2] = 1$ for all $k$.
 The sequence $(a_0, a_1, \ldots)$ has interesting properties: each $a_k$ depends  only on  its successor $a_{k+1}$ and its predecessor $a_{k-1}$ but it is independent of all other variables in the sequence. More precisely,
      $$\mathbb{E}[a_k \ a_j] = 0 \mbox{ if } |k - j| > 1 \mbox{ and } \mathbb{E}[a_k \ a_j] = -1/2 \mbox{ if }|k-j| = 1,\,\, \forall k, j.$$
We now consider the problem of real zeros of the random polynomial $P_n(x) = a_0 + a_1 x + \dots+ a_{n-1} x^{n-1}$. 
From (\ref{ew32wswaws2}) for $H= 0$, 
        $$g(k, k, H) = 1, \,\,  g(k, k+1, H) = g(k+1, k, H) = -1/2  \mbox{ and } g(j, k, H) = 0 \mbox{ for } |k-j| > 1.$$
Then
      \begin{eqnarray*}
 \alpha & = & \sum_{i,j = 0}^{n-1} g(i, j, H)\ x^{i + j} =  \sum_{k=0}^{2n-2}(-1)^k x^k, \\
 \beta & =  & \sum_{i, j=1}^{n-1}  g(i, j, H) \ i \ j \ x^{i+j-2} = \sum_{k=1}^{n-1}k^2 x^{2 k - 2} - 
 \sum_{k=1}^{n-2} k (k + 1) x^{2 k - 1}, \\
 \gamma & = & \sum_{i=0, j=1}^{n-1} g(i, j, H) \ j \ x^{i + j -1} = \frac{1}{2}\sum_{k=1}^{2n-2} (-1)^k  k \ x^{k - 1}. \\
 \end{eqnarray*}
That is, 
     \begin{eqnarray*}
 \alpha & = &\frac{x + x^{2 n}}{x + x^2}, \\
 \beta & =  & -\frac{ x^3 + (n-1) n x^{2 n} + (n-2) n x^{1 + 2 n} - (n-1) n x^{
  2 + 2 n} - (n-1)^2 x^{3 + 2 n}}{(x-1) x^3 (1 + x)^3},\\
 \gamma & = & \frac{-x^2 + (2 n-1) x^{2 n} + 2 (n-1) x^{2 n+1}}{2 x^2 (1 + x)^2}.
\end{eqnarray*}
The average number of real zeros of $P_n(x)$ is 
             $$E_n = \frac{1}{\pi}\int_{-\infty}^\infty \frac{\Delta^{1/2}}{\alpha} dx,\,\,\, \,\,\Delta = \alpha \beta - \gamma^2.$$
We want to show that asymptotically 
         $$E_n \sim (1/\pi) \log n,\,\,\, n\to \infty.$$
First, we consider the number of zeros on the positive real line. Investigating the function $\Delta^{1/2}/\alpha$, yields that for $x\ne \pm 1$, 
         \begin{eqnarray} \label{eqqweq2s3}
      \lim_{n\to \infty} \frac{\Delta^{1/2}}{\alpha} & =&  \frac{\sqrt{(3+x)/(1-x)}}{2+2x} \mbox{ if } |x|<1, \\ \lim_{n\to \infty} \frac{\Delta^{1/2}}{\alpha} & = &  \frac{\sqrt{(1+3x)/(x-1)}}{2x+2x^2} \mbox{ if } |x| >1. 
  \end{eqnarray}
This limit function is integrable on $[0, +\infty)$ but not on $(-\infty, 0]$. In fact
         $$\int_0^{\infty} \lim_{n\to \infty} \frac{\Delta^{1/2}}{\alpha} dx = \pi/3 - \log(2-\sqrt{3}).$$
Moreover denote $f_n(x) = \Delta^{1/2}/\alpha$ and $f(x) = \lim_{n\to \infty} f_n(x)$ for $x\ne 1$. Since $f_n$ and $f$ are all nonnegative, then $|f_n(x) - f(x)| \leq f(x)$ for all $x\ne 1$ and by Fatou's lemma,
 \begin{eqnarray*}
\lim_{n\to \infty} |\int_0^\infty f_n(x) dx - \int_0^\infty f(x) dx| & \leq & \lim_{n\to \infty}  \int_0^\infty |f_n(x) - f(x)| dx\\
 & \leq & \limsup_{n\to \infty} \int_0^\infty |f_n(x) - f(x)| dx\\
 & \leq & \int_0^\infty \limsup_{n\to \infty} |f_n(x) - f(x)| dx\\
 & = & 0.
\end{eqnarray*} 
Therefore
 $$\lim_{n\to \infty} \int_0^\infty f_n(x) dx = \int_0^\infty f(x) dx.$$ That is, 
    $$ \lim_{n\to \infty} \int_0^{\infty} \frac{\Delta^{1/2}}{\alpha} dx = \pi/3 - \log(2-\sqrt{3}).$$
 The interpretation of this limit is that $E_n^+$, the average number of real zeros of the random polynomial $P_n(x)$  in $[0, +\infty)$, is such that $$E_n^+ \sim \pi/3 - \log(2-\sqrt{3}),\,\,n\to \infty.$$ In particular, if we denote by $N_n^+$ the number of positive real zeros of $P_n(x)$, then almost surely, the sequence $(N_1^+, N_2^+, \ldots)$ is bounded.   
This is in sharp contrast to the case where $0 < H < 1$ (in particular to the case where the  random variables $a_i$ are independent for which on average the numbers of positive zeros and negative zeros are equal to $\sim (1/\pi)\log n$ for $n$ large).  \\ 
Moreover, we have the following estimates
\begin{eqnarray*}
\frac{\Delta^{1/2}}{\alpha} &= & \frac{\sqrt{(1+3x)/(x-1)}}{2x+2x^2} + O(n^2 x^{-2n}),\,\,|x|>1\\
\frac{\Delta^{1/2}}{\alpha} &= &  \frac{\sqrt{(3+x)/(1-x)}}{2+2x} + O(n^2 x^{2n}),\,\,|x|<1.
\end{eqnarray*}
For $1< |x| < 1+1/n$, 
\begin{eqnarray*}
\frac{\Delta^{1/2}}{\alpha} =  O(n).
\end{eqnarray*}
First fix $\epsilon >0$, write
  \begin{eqnarray*}
\int_{-\infty}^{-1} \frac{\Delta^{1/2}}{\alpha} dx  = \int_{-\infty}^{-1-n^{\epsilon-1}} \frac{\Delta^{1/2}}{\alpha} dx + \int_{-1 -n^{\epsilon -1}}^{-1-n^{-1}}  \frac{\Delta^{1/2}}{\alpha} dx + \int_{-1-n^{-1}}^{-1} \frac{\Delta^{1/2}}{\alpha} dx.
\end{eqnarray*}
Now
   \begin{eqnarray*}
\int_{-\infty}^{-1-n^{\epsilon-1}} \frac{\Delta^{1/2}}{\alpha} dx & = & \int_{-\infty}^{-1-n^{\epsilon-1}} \frac{\sqrt{(1+3x)/(x-1)}}{2x+2x^2} dx  + O(1)\\
& = & \frac{(1-\epsilon)}{2} \log n + O(1).     
\end{eqnarray*}
Since $\epsilon$ can be taken arbitrarily small, this implies that 
$$ \int_{-\infty}^{-1-1/n} \frac{\Delta^{1/2}}{\alpha} dx = (1/2)\log n + O(1).$$
Also, 
    \begin{eqnarray*}
         \int_{-1+1/n}^0 \frac{\Delta^{1/2}}{\alpha} dx & = & \int_{-\infty}^{-1-1/n}  \frac{\sqrt{(3+x)/(1-x)}}{2+2x}  dx + O(1) = (1/2)\log n + O(1).      
         \end{eqnarray*}
Moreover,
 $$\int_{-1-1/n}^{-1+1/n}\frac{\Delta^{1/2}}{\alpha} dx  = O(1).$$
Finally
 \begin{eqnarray*}
 \int_{-\infty}^{0}\frac{\Delta^{1/2}}{\alpha} dx = (1/2)\log n + O(1).
 \end{eqnarray*}
This concludes the proof of Theorem 2.

\section{Concluding remarks}     
This paper studies random polynomials with correlated random coefficients defined as increments of the classical fractional Brownian motion. These are natural generalisations of random polynomials with independent random coefficients.  We have obtained an asymptotic  formula for the average number of real zeros that generalises  the classical formula of Kac. There are many other properties of zeros of random polynomials only known for the case of independent coefficients. It would be interesting to investigate how these properties can be extended to the case of correlated random variables studied in this paper. One would mention among others the problem of non-zero crossings for random polynomials studied by Dembo and Mukherjee \cite{Dembo_Mukherjee}, the problem of random polynomials having few or no real zeros (Dembo et al. \cite{Dembo_et_al}), the problem of distribution of zeros of random analytic functions (Kabluchko and Zaporozhets \cite{Kabluchko_Zaporozhets}) and the problem of real zeros of linear combinations of orthogonal polynomials (Lubinsky, Pritsker and Xie \cite{Lubinsky}). It would also be interesting to study possible extension the results of Rezakhah and Shemehsavar \cite{Rezakhah1, Rezakhah2} on Brownian motion to the general case of fractional Brownian motion. 

{\bf Acknowledgements.} I thank the anonymous referee whose comments helped to improve the paper. This
research is supported by the Vision Keepers' programme of the University of South Africa.

\end{document}